\documentclass[12pt, a4paper,leqno]{article}
\usepackage[utf8]{inputenc}
\usepackage{amsmath}
\usepackage{amsfonts}
\usepackage{amssymb}
\usepackage{graphicx}
\usepackage{pifont}
\usepackage[margin=3cm]{geometry}
\usepackage{enumerate}
\usepackage {shadow}
\usepackage{amsmath , amsfonts , amssymb , mathrsfs}
\usepackage[linkcolor=blue,colorlinks=true]{hyperref}
\usepackage[T1]{fontenc}
\usepackage{abstract,lipsum}
\usepackage{graphicx}
\usepackage{lmodern}
\usepackage{tocbibind}
\usepackage{microtype}
\usepackage{xcolor}
\usepackage{color}
\usepackage{fancybox}
\usepackage{varioref}
\usepackage{stmaryrd} 
\usepackage{textcomp}
\usepackage{pgf, tikz} 
\usepackage{setspace,booktabs}
\usepackage{enumitem}

\usepackage{xpatch}  
\makeatletter   
\xpatchcmd{\@thm}{\thm@headpunct{.}}{\thm@headpunct{}}{}{}   

\usepackage{fancyhdr}
\numberwithin{equation}{section}
\usepackage{float}
\usepackage[T1]{fontenc}
\usepackage[utf8]{inputenc}

\newtheorem{thm}{Theorem}[section]

\newtheorem{lem}{Lemma}[section]

\newtheorem {prop}{Proposition}[section]
\newenvironment{pr}{{\bf{Proof :~}}} 

\newcommand{\rms}{\mathrm{S}}
\newcommand{\rmi}{\mathrm{I}}

\newcommand{\un}{\underset}
\newcommand{\mbi}{\mathbb{I}}
\newcommand{\mbn}{\mathbf{N}}
\newcommand{\ms}{\mathbf{S}}

\newcommand{\mi}{\mathbf{I}}
\newcommand{\mb}{\mathrm{F}}
\newcommand{\mv}{\mathrm{V}}
\newcommand{\ma}{\mathrm{A}}
\newcommand{\md}{\mathrm{D}}

\usepackage{xcolor}
\usepackage{color}
\usepackage{colortbl}
\newcommand{\n}{\nonumber}

\title{Stochastic and deterministic SIS patch model}
\author{ T. Yeo \thanks{Aix Marseille Univ, CNRS, Centrale Marseille, I2M, Marseille, France; tenan.yeo@univ-amu.fr} \thanks{Univ Félix Houphouët Boigny, Abidjan-Côte d'Ivoire; yeo.tenan@yahoo.fr} }

\begin{document}
	
	\maketitle
	
	\labelformat{rmq}{Remark~#1}
	\labelformat{prop}{Proposition~#1}
	\labelformat{lem}{Lemma~#1}
	\labelformat{thm}{Theorem~#1}
	\labelformat{cor}{Corollary~#1}
	\labelformat{ass}{Assumption~#1}
	\newcommand{\dis}{\displaystyle}
	\newcommand{\tep}{T_{\varepsilon}}
	\newcommand{\sep}{\mathcal{S}_{\varepsilon}}
	\newcommand{\iep}{\mathcal{I}_{\varepsilon}}
	\newcommand{\rep}{\mathcal{R}_{\varepsilon}}
	\newcommand{\T}  {\mathbb{T}}
	\newcommand{\E}  {\mathbb{E}}
	\newcommand{\D}  {\mathbb{D}}
	\renewcommand{\P}  {\mathbb{P}}
	\newcommand{\R}  {\mathbb{R}}
	\newcommand{\bV}  {\big\Vert}
	\newcommand{\BV}  {\Big\Vert}
	\newcommand{\sbs}{\subset}
	\newcommand{\mfrm}{\mathfrak{M}}
	\newcommand{\nr}{\mathrm{N}}
	\newcommand{\nn}  {\noindent}
	\newcommand{\fpr} {\begin{flushright}
	$\square$
	\end{flushright}}
\newcommand{\mtt}{\mathtt{T}}
	\newcommand{\fprb} {\begin{flushright}
		$\square$
\end{flushright}}
	
\begin{abstract}
\nn  Here, we consider  a SIS epidemic model where the individuals are distributed on several distinct patches. We construct  a stochastic  model and then prove that it converges to a deterministic model as the total population size tends to infinity. Furthermore we show the existence and the global stability of a unique endemic equilibrium  provided that the migration rates of susceptible and infectious individuals are equal. Finally we compare the equilibra  with those of the homogeneous model, and with  those  of isolated patches.
\end{abstract}

\setcounter{section}{-1}

\section{Introduction}	

Early epidemic models were formulated assuming that individuals in the population mix homogeneously
\cite{ab2000,nj,BP2019,wo,gh}. In this consideration, all pairs of individuals in the population have the same probability of coming into contact with each other. But, it is well known that in a large population several groups can be formed due to heterogeneity arising, for example, from social and economic factors. Some people may live in cities while others may live in rural areas. Consequently, demographic and disease parameters may vary for each group, and then the persistence and extinction  of infectious diseases in those communities can be different.  Furthermore, people may travel between the groups, which leads to the spread of the disease between groups. Then it is clear that spatial  heterogeneity, habitat connectivity and rates of movement of individuals play an important role in the outbreak of an infectious disease. Those considerations have lead to the development of epidemic models that take into account the structure of the population. For this reason, several authors studied heterogeneous epidemic models, by structuring the spatial environment into patches. Fulford et al. \cite{fr} developed a Susceptible-Exposed-Infectious (SEI) metapopulation model for the spread of an infectious agent by migration.  Salmani $\&$ Van den Driessche \cite{sl} considered a SEIRS deterministic model in which travel rates were assumed to depend on the disease status. Disease spread in metapopulation models involving discrete patches has been also investigated by  Arino et al. \cite{al}, Arino $\&$ Van den Driessche \cite{ar},  Jin $\&$ Wang \cite{jw}, Wang $\&$ Mulone \cite{wm}, Wang $\&$ Zho \cite{wz}. Arino $\&$ Van den Driessche  developed a general framework for movement of susceptible, exposed, infectious, and recovered individuals (SEIRS model) and define a mobility matrix, an irreducible matrix that defines the spatial arrangement of patches and rates of movement between patches. Wang and colleagues studied uniform persistence and global stability of disease-free and endemic equilibria in Susceptible-Infectious-Susceptible (SIS) metapopulation models. In a similar setting, Allen et al. \cite{al7} showed for a  SIS deterministic patch model that, while the population is at an endemic level and if infectious individuals travel between the patches but the rate of travel for susceptible individuals approaches zero, then,  the endemic equilibrium approaches a disease–free equilibrium. \\
\nn  Deterministic patch models describe the spread of disease under the assumptions of mass action, relying on the law of large numbers. But the most natural way to describe the spread of disease is stochastic.
 The probabilistic point of view is recent. Let us mention some authors which treated stochastic epidemic models. In 2007, McCormack $\&$ Allen \cite{mck} studied a SIR and a SIS epidemic model. In this setting  both models are deterministic and stochastic. They showed that travel between patches can lead to either disease persistence or extinction in all patches. Considering stochastic model, Clancy \cite{cl} proposed a SIR model, and then showed that movement of infectious individuals can decrease the spread of the disease. In the same considerations, Sani et al. \cite{san} introduced a multi-group SIR model for the spread of AIDS. In this study, the authors used Markov process to describe the model. Using an approximating system of the ODEs, they analysed the equilibrium behaviour of the stochastic model. Finally, let's mention that stochastic epidemics in a homogeneous community has been studied recently by Britton $\&$ Pardoux \cite{BP2019}.  
 
 \medskip
 
The rest of the paper is organised as follows. In section 1 we introduce  a stochastic  model on a finite number of patches. Section 2 is devoted to the law of large numbers. In section 3 we show that the limit deterministic model has an unique  endemic equilibrium (EE) which is gloably asymptotically stable.  Finally we compare this endemic equilibrium with the one of the homogeneous model in the section 4.
\section{The model}
Consider a population consisting of N individuals, where each individual is located at one of $j$ geographically distinct patches. Sites (or vertices) represent human communities in which the disease can diffuse and grow. The edges represent links between communities (see figure 1 below). Individuals in that population can be  classified according to their ability to transmit the disease to others. Susceptible individuals are those who do not have the disease and  who can become infected. Infectious individuals are those  who are infected by the disease and can transmit it to susceptible individuals. In this work, attention is given to the SIS model, but the same approach can be developed in the case of the SIRS model and of the SIR model with demography. For any patch $j$, the transmission of the disease depends on three factors: the rate of contacts, the probability that a contact is made with a susceptible individual, and the probability that a contact between an infectious and a susceptible individual leads to a successful transmission (see e.g \cite{mb,BP2019}).
$\rms_j(t)$ (resp. $\rmi_j(t)$ ) denotes the number of susceptible (resp. infectious)  individuals in patch $j$ at time $t.$ We formulate a random Markov epidemic model as a Poisson process driven stochastic differential equation (SDE). In what follows the $\mathrm{P}_{\!\!j}$ are mutually independent standard Poisson processes. Infections, healings and migrations of individuals happen according to  Poisson processes.
In this model
\begin{enumerate}
	\item[$\bullet$]infections are local;
	\item[$\bullet$]when an infectious individual cures, he immediately becomes susceptible again;
	\item[$\bullet $] each  infectious individual meets other individuals at some rate $\alpha_j$. The encounter results in a new infection with probability $p_j$ if the patner of the encounter is susceptible, which happens with probability $\dis \dfrac{\rms_j(t)}{\rms_j(t)+\rmi_j(t)}$, since we assume that individuals in each patch mix homogeneously. Letting $ \lambda_j= \alpha_jp_j$  and summing over the infectious individuals at time $t$ gives the rate $\dis \lambda_j \dfrac{ \rms_j(t)}{\rms_j(t)+\rmi_j(t)}\rmi_j(t)$ at time $t$.  Then $\dis \mathrm{P}_{\!\!j}^{inf}\left(  \int_0^t \lambda_j \dfrac{ \rms_j(r)\rmi_j(r)}{\rms_j(r)+\rmi_j(r)}dr \right)$ counts the number of transitions of type $\rms\longrightarrow \rmi$ on the patch $j$ between time $0$ and time $t$;
	\item[$\bullet $] recovery of an infectious happens at rate $ \gamma_j $, so $ \dis \mathrm{P}_{\!\!j}^{rec}\left(\int_0^t \gamma_j \rmi_j(r)dr \right)$ counts the number of  transitions of type  $ \rmi \longrightarrow \rms $ on the patch $j$ between time $0$ and time~$ t$.  
	\item[$\bullet $] The term $ \dis \mathrm{P}_{\!\!S,j, k}^{mig}\left(\int_0^t \nu_S a_{jk} \rms_j(r)dr \right) $ counts  the number of migrations of susceptible individuals from patch $j$ to  $k$, if we assume that each susceptible migrates  from $j$ to $k$  at rate $ \nu_Sa_{jk}$, and similarly for the compartment $\rmi$.
\end{enumerate}

\vspace{-0.2cm}

Here, $\nu_S$ and $ \nu_I$ are the diffusion coefficients for susceptible and infectious individuals, respectively. $a_{ij}$ represents the degree of movement from patch $i$ into patch $j$.

\newcommand{\cercle}[2]{\draw[thick,green,fill=white] (#1,#2) circle(.6);}
\newcommand{\cercleb}[2]{\fill[blue] (#1,#2) circle(.1);}
\newcommand{\cercler}[2]{\fill[red] (#1,#2) circle(.1);}
\begin{center}
\begin{tikzpicture}
\draw[thick,blue] (0,0)--(4.2,0.5)--(2,1.5)--(0,0)--(0,2.5)--(1.5,4.5)--(4.2,0.5);
\draw[thick,blue] (0,0)--(1.5,4.5)--(2,1.5);
\draw[thick,blue] (2,1.5)--(4.2,2.5);
\draw[thick,blue] (4.2,2.5)--(1.5,4.5);
\draw[thick,blue] (4.2,2.5)--(5.5,4.5);
\draw[thick,blue] (5.5,4.5)--(2,1.5);
\draw[thick,blue] (5.5,4.5)--(1.5,4.5);
\draw[thick,blue] (5.5,4.5)--(0,2.5);
\draw[thick,blue] (5.5,4.5)--(9,2);
\draw[thick,blue] (4.2,0.5)--(9,2);
\draw[thick,blue] (0,2.5)--(9,2);
\draw[thick,blue] (7,-2)--(9,2); 
\draw[thick,blue] (7,-2)--(4.2,0.5);
\draw[thick,blue] (7,-2)--(5.5,4.5);
\draw[thick,blue] (7,-2)--(2.5,-3);
\draw[thick,blue] (2.5,-3)--(0,0);
\draw[thick,blue] (2.5,-3)--(4.2,0.5);
\draw[thick,blue] (0,0)--(7,-2);
\draw[thick,blue] (2.5,-3)--(9,2);
\draw[thick,blue] (8.5,5)--(9,2);
\draw[thick,blue] (8.5,5)--(4.2,0.5);
\draw[thick,blue] (8.5,5)--(4.2,2.5);
\draw[thick,blue] (8.5,5)--(5.5,4.5);
\draw[thick,blue] (4.5,7)--(5.5,4.5);
\draw[thick,blue] (4.5,7)--(8.5,5);
\draw[thick,blue] (4.5,7)--(9,2);
\draw[thick,blue] (4.5,7)--(4.2,0.5);
\draw[thick,blue] (4.5,7)--(1.5,4.5);
\draw[thick,blue] (2.2,-0.5)--(2.5,-3);
\draw[thick,blue] (2.2,-0.5)--(4.2,0.5);
\draw[thick,blue] (2.2,-0.5)--(0,2.5);
\draw[thick,blue] (2.2,-0.5)--(2,1.5);
\draw[thick,blue] (3,5.7)--(2,1.5);
\draw[thick,blue] (3,5.7)--(8.5,5);
\draw[thick,blue] (3,5.7)--(5.5,4.5);
\draw[thick,blue] (3,5.7)--(4.2,2.5);
\draw[thick,blue] (7.95,0)--(4.2,0.5);
\draw[thick,blue] (7.95,0)--(3,5.7);
\draw[thick,blue] (7.95,0)--(8.5,5);
\draw[thick,blue] (7.95,0)--(4.2,0.5);
\draw[thick,blue] (7.95,0)--(4.5,-1.3);
\cercle{0}{0};
\cercle{0}{2.5};

\cercle{3}{5.7}
\cercleb{3}{5.7}
\foreach \k in {0,...,1}{
	\cercler{{3+.4*cos(\k*120)}}{{5.7+.4*sin(\k*120)}}
}
\foreach \k in {2}{
	\cercleb{{3+.4*cos(\k*120)}}{{5.7+.4*sin(\k*120)}}
}
\foreach \k in {0,...,1}{ \cercleb{{3+.4*cos(\k*120+60)}}{{5.7+.4*sin(\k*120+60)}} }
\foreach \k in {4}{ \cercleb{{3+.4*cos(\k*120+60)}}{{5.7+.4*sin(\k*120+60)}} }
\foreach \k in {5}{ \cercleb{{3+.2*cos(\k*55+70)}}{{5.7+.3*sin(\k*20+160)}} }

\cercle{6}{2.5}
\cercler{6}{2.5}
\foreach \k in {0,...,2}{
	\cercler{{6+.4*cos(\k*120)}}{{2.5+.4*sin(\k*120)}}
}
\foreach \k in {0,...,1}{ \cercleb{{6+.4*cos(\k*120+60)}}{{2.5+.4*sin(\k*120+60)}} }
\foreach \k in {4}{ \cercleb{{6+.4*cos(\k*120+60)}}{{2.5+.4*sin(\k*120+60)}} }
\foreach \k in {5}{ \cercleb{{6+.2*cos(\k*55+70)}}{{2.6+.3*sin(\k*20+160)}} }

\cercle{2.2}{-0.5}
\cercleb{2.2}{-0.5}
\foreach \k in {0,...,2}{
	\cercleb{{2.2+.4*cos(\k*120)}}{{-0.5+.4*sin(\k*120)}}
}
\foreach \k in {0,...,1}{ \cercleb{{2.2+.4*cos(\k*120+60)}}{{-0.5+.4*sin(\k*120+60)}} }
\foreach \k in {4}{ \cercleb{{2.2+.4*cos(\k*120+60)}}{{-0.5+.4*sin(\k*120+60)}} }
\foreach \k in {5}{ \cercleb{{2.2+.2*cos(\k*55+70)}}{{-0.5+.3*sin(\k*20+160)}} }

\cercle{4.5}{-1.3}
\cercleb{4.5}{-1.3}
\foreach \k in {0,...,1}{
	\cercler{{4.5+.4*cos(\k*120)}}{{-1.3+.4*sin(\k*120)}}
}
\foreach \k in {2}{
	\cercleb{{4.5+.4*cos(\k*120)}}{{-1.3+.4*sin(\k*120)}}
}
\foreach \k in {0,...,1}{ \cercleb{{4.5+.4*cos(\k*120+60)}}{{-1.3+.4*sin(\k*120+60)}} }
\foreach \k in {4}{ \cercleb{{4.5+.4*cos(\k*120+60)}}{{-1.3+.4*sin(\k*120+60)}} }
\foreach \k in {5}{ \cercleb{{4.5+.2*cos(\k*55+70)}}{{-1.3+.3*sin(\k*20+160)}} }

\cercle{6.4}{0.2}
\cercler{6.4}{0.2}
\foreach \k in {0,...,5}{
	\cercler{{6.4+.4*cos(\k*120)}}{{0.2+.4*sin(\k*120)}}
}

\cercleb{0}{0}
\foreach \k in {0,...,2}{
	\cercler{{0+.4*cos(\k*120)}}{{.4*sin(\k*120)}}
}
\foreach \k in {0,...,1}{ \cercleb{{.4*cos(\k*120+60)}}{{.4*sin(\k*120+60)}} }
\foreach \k in {0,...,1}{ \cercleb{{.4*cos(\k*120+60)}}{{.4*sin(\k*120+60)}} }
\foreach \k in {4}{ \cercleb{{0+.4*cos(\k*55+70)}}{{.4*sin(\k*20+170)}} }

\cercle{0}{2.5}
\cercler{0}{2.5}
\foreach \k in {0,...,2}{
	\cercler{{0+.4*cos(\k*120)}}{{2.5+.4*sin(\k*120)}}
}
\foreach \k in {0,...,1}{ \cercleb{{+.4*cos(\k*120+60)}}{{2.5+.4*sin(\k*120+60)}} }
\foreach \k in {4}{ \cercleb{{0+.4*cos(\k*55+70)}}{{2.5+.4*sin(\k*20+170)}} }

\cercle{1.5}{4.5}
\cercler{1.5}{4.5}
\foreach \k in {0,...,4}{
	\cercleb{{1.5+.4*cos(\k*120)}}{{4.5+.4*sin(\k*120)}}
}
\foreach \k in {0,...,1}{ \cercleb{{1.5+.4*cos(\k*120+60)}}{{4.5+.4*sin(\k*120+60)}} }
\foreach \k in {4}{ \cercleb{{1.5+.4*cos(\k*55+70)}}{{4.5+.4*sin(\k*20+170)}} }

\cercle{5.5}{4.5}
\cercleb{5.5}{4.5}
\foreach \k in {0,...,4}{
	\cercleb{{5.5+.4*cos(\k*120)}}{{4.5+.4*sin(\k*120)}}
}
\foreach \k in {0,...,1}{ \cercleb{{5.5+.4*cos(\k*120+60)}}{{4.5+.4*sin(\k*120+60)}} }
\foreach \k in {4}{ \cercleb{{5.5+.4*cos(\k*55+70)}}{{4.5+.4*sin(\k*20+170)}} }

\cercle{9}{2}
\cercleb{9}{2}
\foreach \k in {0,...,4}{
	\cercleb{{9+.4*cos(\k*120)}}{{2+.4*sin(\k*120)}}
}
\foreach \k in {0}{ \cercler{{9+.4*cos(\k*120+60)}}{{2+.4*sin(\k*120+60)}} }
\foreach \k in {5}{ \cercleb{{9+.4*cos(\k*45+70)}}{{2+.4*sin(\k*20+170)}} }

\cercle{7}{-2}
\cercleb{7}{-2}
\foreach \k in {0,...,4}{
	\cercleb{{7+.4*cos(\k*120)}}{{-2+.4*sin(\k*120)}}
}
\foreach \k in {0,...,1}{ \cercleb{{7+.4*cos(\k*120+60)}}{{-2+.4*sin(\k*120+60)}} }
\foreach \k in {4}{ \cercleb{{7+.4*cos(\k*55+70)}}{{-2+.4*sin(\k*20+170)}} }

\cercle{2.5}{-3}
\cercleb{2.5}{-3}
\foreach \k in {0,...,3}{
	\cercleb{{2.5+.4*cos(\k*120)}}{{-3+.4*sin(\k*120)}}
}
\foreach \k in {0,...,1}{ \cercler{{2.5+.4*cos(\k*120+60)}}{{-3+.4*sin(\k*120+60)}} }
\foreach \k in {4}{ \cercleb{{2.5+.4*cos(\k*55+70)}}{{-3+.4*sin(\k*20+170)}} }

\cercle{8.5}{5}
\cercler{8.5}{5}
\foreach \k in {0,...,3}{
	\cercleb{{8.5+.4*cos(\k*120)}}{{5+.4*sin(\k*120)}}
	\cercler{{8.5+.4*cos(\k*120+60)}}{{5+.4*sin(\k*120+60)}} }

\cercle{4.5}{7}
\foreach \k in {0,...,5} {
	\cercleb{{4.5+.4*cos(\k*60)}}{{7+.4*sin(\k*60)}}
}
\foreach \k in {0,1,2} {
	\cercler{{4.5+.2*cos(\k*120+90)}}{{7+.2*sin(\k*120+90)}}
}

\cercle{4.2}{0.5}
\foreach \k in {0,...,4} {
	\cercleb{{4.2+.4*cos(\k*72)}}{{.5+.4*sin(\k*72)}}
}
\foreach \k in {0,...,3} {
	\cercler{{4.2+.18*cos(\k*90+18)}}{{.5+.18*sin(\k*90+18)}}
}

\cercle{4.2}{2.5}
\cercler{4.2}{2.5}
\foreach \k in {0,...,3} {
	\cercleb{{4.2+.4*cos(\k*120)}}{{2.5+.4*sin(\k*120)}}
	\cercler{{4.2+.4*cos(\k*120+60)}}{{2.5+.4*sin(\k*120+60)}}
}

\cercle{2}{1.5}
\cercleb{2}{1.5}
\foreach \k in {0,...,6} {
	\cercleb{{2+.4*cos(\k*60)}}{{1.5+.4*sin(\k*60)}}
}

\cercle{7.95}{0}
\cercler{7.95}{0}
\foreach \k in {0,...,3} {
	\cercleb{{7.95+.4*cos(\k*120)}}{{.4*sin(\k*120)}}
	\cercleb{{7.95+.4*cos(\k*120+60)}}{{.4*sin(\k*120+60)}}
	
}
\end{tikzpicture}
\hspace{-2.55cm}{\color{red}$\bullet$}  infectious individuals

\hspace{8cm}{\color{blue}$\bullet$} susceptible individuals\\

\hspace{-4cm}Figure 1: Metapopulation
\end{center}
 Then the  propagation of the illness can be modeled by  the following  system of stochastic differential equations
\begin{equation}\label{sm}
\hspace{-0.5cm}\left\{ 
\begin{aligned} 
\bigskip 
\rms_j(t) &=  \rms_j(0) - \mathrm{P}_{\!\!j}^{inf}\left(  \int_0^t \lambda_j \dfrac{ \rms_j(r)\rmi_j(r)}{\rms_j(r)+\rmi_j(r)}dr \right) +\mathrm{P}_{\!\!j}^{rec}\left(\int_0^t \gamma_j \rmi_j(r)dr \right) \\
& - \sum_{\un{k\ne j}{k=1}}^{\ell}\mathrm{P}_{\!\!S,j,k}^{mig}\left(\int_0^t \nu_S a_{jk} \rms_j(r)dr \right) + \sum_{\un{k\ne j}{k=1}}^{\ell}\mathrm{P}_{\!\!S,k,j}^{mig}\left( \int_0^t \nu_S a_{kj} \rms_k(r)dr \right)
\\[2mm]
\rmi_j(t) & =  \rmi_j(0) +\mathrm{P}_{\!\!j}^{inf}\left(\int_0^t \lambda_j \dfrac{ \rms_j(r)\rmi_j(r)}{\rms_j(r)+\rmi_j(r)}dr \right)-\mathrm{P}_{\!\!j}^{rec}\left(\int_0^t \gamma_j \rmi_j(r)dr \right)\\ 
& - \sum_{\un{k\ne j}{k=1}}^{\ell}\mathrm{P}_{\!\!I,j,k}^{mig}\left(\int_0^t \nu_I a_{jk} \rmi_j(r)dr \right) + \sum_{\un{k\ne j}{k=1}}^{\ell}\mathrm{P}_{\!\!I,k,j}^{mig}\left( \int_0^t \nu_I a_{kj} \rmi_k(r)dr \right)\\
t &\in  [0,T], \; \; \; j=1,\cdots, \ell.
\end{aligned}
\right. 
\end{equation}

 In the next section we show that this stochastic model converges to a deterministic epidemic patch model as the total size of population tends to infinity.
\section{Law of large numbers} 
We introduce the martingales $\mathrm{M}_j(t)=\mathrm{P}_{\!\!j}(t) - t $  and we  look instead at the renormalized model by dividing the size of the population in each compartment by $ \nr $. Hence by setting \\

$ \rms_j^{\nr}(t) = \dfrac{\rms_j(t)}{\nr}$,  \; \;   $ \rmi_j^{\nr}(t) = \dfrac{\rmi_j(t)}{\nr}$, \; \;
$\dis \rms^\nr(t) = \left( \begin{array}{cl}
\rms_1^\nr(t)\\
\vdots \\
\rms_{\ell}^\nr(t)
\end{array}
\right)$,  \; \;  $ \rmi^\nr(t) = \left( \begin{array}{cl}
\rmi_1^\nr(t)\\
\vdots \\
\rmi_{\ell}^\nr(t) 
\end{array} 
\right)  $, and \\  $\dis  \mathcal{Z}^\nr(t)= \left( \begin{array}{cl}
\rms^\nr(t)\\
\rmi^\nr(t)
\end{array}
\right),$ then the  stochastic model takes the aggregated form 
\begin{eqnarray}\label{fag}
\mathcal{Z}^\nr(t) 
&=& \mathcal{Z}^\nr(0)+ \int_0^t b\big(r,\mathcal{Z}^\nr(r)\big)dr + \sum_{j=1}^k\frac{h_j}{\nr} \mathrm{M}_j\Bigg( \nr\int_0^t \beta_j\big( \mathcal{Z}^\nr(r)\big)dr\Bigg),  
\end{eqnarray}
where  $ k$ is the total number of $\mathrm{P}_{\!\!j}$'s in the system, and     
\begin{eqnarray}\label{b}
 b\big(r,\mathcal{Z}^\nr(r)\big)=  \sum_{j=1}^kh_j \beta_j\big( \mathcal{Z}^\nr(r)\big);   
 \end{eqnarray}
the vectors $h_j \in \{ -1 , 0 , 1\}^{2\ell}$ denote the respective jump directions with jump rates $\beta_j$ . The rates 
\begin{eqnarray}\label{tau}
 \beta_.\big(\mathcal{Z}^\nr(t)\big)\!\!\!& \in &\!\!\!\Big\{\; \dfrac{ \lambda_j\rms_j^\nr(t)\rmi_j^\nr(t)}{\rms_j^\nr(t)+\rmi_j^\nr(t)} , \;   \gamma_j \rmi_j^\nr(t),  \; \nu_S\, a_{ij}\rms_j^\nr(t), \\
 &&  \; \; \; \; \; \; \; \; \;\; \; \; \; \; \;\; \; \; \; \; \; \; \; \; \nu_I\, a_{ij} \rmi_j^\nr(t),  \; i,j\in\{ 1,\cdots, \ell\}\; \Big\}.\nonumber
 \end{eqnarray}
 Concerning the initial condition, we assume that $\mathcal{Z}^\nr(0)=z_\nr=[\nr x]/\nr$, for some $x\in [0,1]^{\ell}$, where $[\nr x]$ is a vector of integers.
 
 We set  $ \dis  \mathcal{F}_t^\nr = \sigma\{\; \mathcal{Z}_j^\nr(r) , \;  0\leq r \leq t,\;  j=1,\cdots, \ell \; \} $ and  we shall assume that the process $\{\mathcal{Z}^\nr(t), t\ge 0\}$ is defined on the filtered probability space $\dis \big(\Omega, \mathcal{F}, \mathcal{F}_t^\nr, \P\big)$. In what follows, $\Vert u \Vert $ denotes the $L^1$ norm  of an $\ell$–dimensional vector u. More precisely, $\dis \Vert u \Vert = \sum_{j=1}^{\ell} \vert u_j \vert $. We shall say that a vector $u$ is nonnegative (resp. positive) if all its elements are nonnegative (resp. positive), in with case  we will write $ u\ge 0$ (resp. $u> 0$).  The following theorem shows that the solution of the stochastic model $(\ref{fag})$ converges a.s. locally
uniformly in $t$ to the solution of a deterministic model, as the total population size  $\nr$ tends to infinity. 
\begin{thm} \label{llns} 
	$\bf{[Law \ of \ Large \ Numbers]}$ \\ ~
	Let  $\mathcal{Z}^\nr $ denote the  solution  of the  SDEs (\ref{fag})  and $ z $ the unique solution of the system of ordinary differential equations \, $\dfrac{dz}{dt}(t)=  b(t,z(t))$, $z(0)=x$. \\
	Let us fix an arbitrary $T > 0$.
	Then $ \dis   \underset{0\leq t\leq T}{\sup}\Big\Vert  \mathcal{Z}^\nr(t)-z(t) \Big\Vert \longrightarrow 0  \; \text{a.s.} ,  \; as \; \;   \nr \rightarrow + \infty .$ 	
\end{thm}
Note that the solution $z(t)=\big(\ms_1(t),\mi_1(t), \ms_2(t),\mi_2(t),\cdots, \ms_{\ell}(t),\mi_{\ell}(t)\big)^{\!\mtt}$ of the deterministic model satisfy
 \begin{equation}\label{md}
 \left\{ 
 \begin{aligned} 
 \frac{d\,\ms_j}{dt}(t) & \!=\!  - \lambda_j\dfrac{\ms_j(t)\mi_j(t)}{\ms_j(t)+\mi_j(t)} +\gamma_j\,\mi_j(t)+ \nu_S \sum_{\un{k\ne j}{k=1}}^{\ell}\Big(a_{kj}\,\ms_k(t)-a_{jk}\,\ms_j(t) \Big) \\ 
 \frac{d\,\mi_j}{dt}(t) & \! =\!  \lambda_j\dfrac{\ms_j(t)\mi_j(t)}{\ms_j(t)+\mi_j(t)} -\gamma_j\,\mi_j(t) +\nu_I \sum_{\un{k\ne j}{k=1}}^{\ell}\Big(a_{kj}\,\mi_k(t)-a_{jk}\,\mi_j(t) \Big)\\
 \ms_j(0)&\ge 0 , \; \mi_j(0) \ge 0 \\
 j&= 1, \cdots , \ell .
 \end{aligned} \right.
 \end{equation}
 $\ms_j(t)$ (resp. $\mi_j(t)$ ) is the proportion of the  total susceptible (resp. infectious)  population which is localized  on the site $ j$ at time $ t$.

\ref{llns} ensures that, as the population size $ \nr $ becomes large, the proportion of susceptible and infectious individuals at each patch is well approximated, on any bounded intervall $[0,T]$ by the solution of the ODEs (\ref{md}),  provided the scaled process starts close to an initial value of the ODEs.

\bigskip

\ref{llns} is a special case of  Theorem 2.2.7 of Britton $\&$ Pardoux \cite{BP2019}, where the proof written for the homogeneous model covers our situation as well. One of the earliest references on this convergence result is Ethier $\&$ Kurtz \cite{ku} (chapter 11, Theorem 2.1).  Thus, we do not give details and  refer the reader to those papers for a complete proof. We briefly sketch the idea of the proof. First note that $ 0 \le \mathcal{Z}^\nr(t)\le 1$, for all $t\in [0,T]$. By using the law of large numbers for Poisson processes and the second Dini Theorem,  it follows that 
$$   \underset{0\leq t\leq T}{\sup}{\Big\rvert \sum_{j=1}^{\ell}\dfrac{h_j}{\nr}\mathrm{M}_j\left(\nr\int_0^t \beta_j(\mathcal{Z}^\nr(r))dr\right)\Big\rvert}\overset{a.s.}{\longrightarrow}0 .$$
Next, it is not hard to see that  $ b(t,z)$ is a globally Lipchitz function  of $z$, locally uniformly in $t$.  From this fact, it follows that, for all $ t\in[0,T]$, 
\begin{eqnarray}\label{gl}
\Big\Vert \mathcal{Z}^\nr(t)- z(t)\Big\Vert &\le & \Big\Vert z_\nr- x \Big\Vert + \Big\rVert \sum_{j=1}^{\ell}\dfrac{h_j}{\nr}\mathrm{M}_j\left(\nr\int_0^t \beta_j(\mathcal{Z}^\nr(r))dr\right)\Big\rVert \\
&+&  C\!\!\int_0^t \Big\Vert \mathcal{Z}^\nr(r)- z(r)\Big\Vert dr, \n
\end{eqnarray} 
where $C$ is the Lipschitz constant  of $b$.
Finally, the result follows from  Gronwall's Lemma and the fact that the two first terms  in  the right-hand side of (\ref{gl}) tend to zero as $\nr \to \infty$.

\smallskip

In the following section, we study the equilibra of this system of  ODEs.
\section{Equilibra  of the ODEs and their stability}
In this section, we consider properties of the disease free equilibrium (DFE) and the endemic equilibrium (EE), including its existence, uniqueness and stability. Throughout this section, we assume that the connectivity matrix $\mathbb{A}=\big(a_{ij}\big)_{\un{1\le j\le \ell}{1\le i\le \ell}}$ is irreducible  and symmetric. This irreducibility assumption implies that the patches cannot be separated into two disjoint subsets such that there is no migration of individuals from one subset to the other. That is, for any $j,k \in \{ 1, \cdots,\ell\}$, $j\ne k$, there exists $s\ge 2$, a sequence $ j_1, j_2,\cdots, j_s \in \{ 1, \cdots,\ell\}$ such that $j_1=j$,  $j_s=k$ and $a_{j_ij_{i+1}}\ne 0$, $\forall i\in \{1, \cdots,s-1\}$. We shall  say that a matrix $M=\big(m_{ij}\big)_{\un{1\le j\le \ell}{1\le i\le \ell}}$ is nonnegative (resp. positive) if all its elements are nonnegative (resp. positive), in with case  we will write $ M\ge 0$ (resp. $M> 0$). In what follows, we  set  $\mbn_j(t) =\ms_j(t) + \mi_j(t) $.
 
Let us mention that the system of ODEs (\ref{md}) obtained in \ref{llns} has been studied by Allen et al. \cite{al7},  where the authors studied  the asymptotic profiles of the steady states. First, using the irreducibility of the connectivity  matrix, they show
\begin{lem}[Allen \cite{al7}\,]$\bf{[Existence\ and \ uniqueness \ of \ the \ DFE]}$  ~\\
	The system (\ref{md}) has a unique  disease-free equilibrium which is  given by 
	$$ \hat{z}:=\Big(\hat{\ms}_1, \hat{\mi}_1,\hat{\ms}_2, \hat{\mi}_2, \cdots, \hat{\ms}_{\ell}, \hat{\mi}_{\ell}\Big)=\Big(\dfrac{1}{\ell}, 0, \dfrac{1}{\ell}, 0, \cdots, \dfrac{1}{\ell}, 0 \Big).$$
\end{lem}
The DFE always exists, an important question is whether an outbreak of the disease can occur when the population initially contains a small number of infected individuals.
This question may be addressed using stability analysis of the DFE. In fact if $\mathcal{R}_0$, the basic reproduction number (the expected number of secondary cases produced, in a completely susceptible population, by a typical infective individual) is less than 1, the DFE is globally asymptotically stable. That is,  the trajectory of the ODEs which starts close to the DFE will be attracted towards the DFE.  That is the content of the next lemma. Before, we use the next-generation matrix approach of Van den Driessche $\&$ Watmough \cite{vd} to compute $\mathcal{R}_0$.  Define the matrices
$$ \ma=\text{diag}(\gamma_j)_{1\le j\le\ell}, \;  \; \; \mb= \text{diag}(\lambda_j)_{1\le j\le\ell}\; \text{and}\;  \; \md=(d_{ij})_{1\le i,\, j\le \ell}$$ with 
$\dis d_{ij}=\left\{ \begin{array}{lr}
	\dis -\sum_{\un{k\ne i}{k=1}}^{\ell} a_{ik} & \; \; \text{if} \; \; \;  i=j, \\
	\; \; \; \; a_{ij} & \; \; \text{if} \; \; \; i\ne j. 
\end{array}\right.$ 

\smallskip
A direct application of the result of the above reference yields the following Proposition.
\begin{prop}
The basic reproduction number for (\ref{md}) is the spectral radius of the next-generation matrix: 
$$ \mathcal{R}_0 = \rho(-\mb\mv^{-1}),$$
where   $\mv= \nu_I \md-\ma$.
\end{prop} 

We have also the 
\begin{lem}[Allen et al. \cite{al7}\,]$\bf{[Stability\  of \ the \ DFE]}$~\\
	If $\mathcal{R}_0 < 1$, then the disease-free equilibrium $\hat{z}$ is globally asymptotically stable, that is  $$z(t)\longrightarrow \hat{z},\; \; \text{ as} \; \;  t\to\infty.$$
\end{lem}

We will show at the end of this section that, under a specific condition, the disease-free equilibrium is globally asymptotically stable if $ \mathcal{R}_0=1$. 
The existence and uniqueness of the EE for (\ref{md}) is shown in \cite{al7} under the assumption that $\mathcal{R}_0>1$. In that work, the authors were not able to prove the stability of the EE, but conjectured that this EE attracts all solutions whose initial conditions have a nonzero proportion of infectious (and numerical simulations suggest that this is indeed the case).  Here, we employ the approach in Bichara et al. \cite{ba} to prove the globally stability of the EE, under the assumption  that  infectious and susceptible individuals have the same diffusion rate $ \nu_{S}=\nu_{I}:=\nu$.\\ Assuming that $ \nu_{S}=\nu_{I}:=\nu$, then the sytem given by (\ref{md}) is equivalent to
\begin{equation}\label{detn}
\hspace{-0.2cm}\left\{ 
\begin{aligned} 
\frac{d\,\mbn_j(t)}{dt}&=\nu\sum_{\un{k\ne j}{k=1}}^{\ell}\Big(a_{kj}\mbn_k(t)-a_{jk}\mbn_j(t) \Big) \\
\frac{d\,\mi_j(t)}{dt}&=\lambda_j\Big(  1 -\frac{\mi_j(t)}{\mbn_j(t)}\Big)\mi_j(t)-\gamma_j\,\mi_j(t) 
+\nu\sum_{\un{k\ne j}{k=1}}^{\ell}\Big(a_{kj}\mi_k(t)-a_{jk}\mi_j(t) \Big) \\
j&= 1, \cdots , \ell,
\end{aligned} 
\right.
\end{equation}
which can be  written in the form 
\begin{equation}\label{deta}
\hspace{-0.2cm}\left\{ 
\begin{aligned} 
\frac{d\,\mbn(t)}{dt}& =\nu \md\mbn(t)\\
\dfrac{d\,\mi(t)}{dt}& = \nu \md\,\mi(t)- \ma\,\mi(t)+ \Big(\mbi_\ell-\text{diag}(\mbn_j^{-1}(t))\text{diag}(\mi(t))\Big)\mb\,\mi(t),
\end{aligned} 
\right.
\end{equation}
where $ \mbn= (\mbn_1, \cdots, \mbn_\ell)^{\mtt}$, $\mi=(\mi_1, \cdots, \mi_\ell)^{\mtt}$ and $ \mbi_\ell$ is the identity matrix with dimension $\ell \times \ell .$ Note that $ \Big(\mbi_{\ell}-\text{diag}(\mbn_j^{-1}(t))\text{diag}(\mi(t))\Big)\mb\,\mi(t) $ is the vector of new infections, $  \nu \md\mi(t)$ is the vector of migrations of infectious individuals and $ \ma\mi(t) $ is the vector of transitions of individuals from the compartment $ \mi $ to  the compartment $\ms$.

\begin{lem}\label{stn}
	The system $ \dis \frac{d\,\mbn(t)}{dt}=\nu \md\mbn(t)$ has a unique global asymptotically stable equilibrium.
\end{lem}

\begin{pr} Let $Q$ be the matrix such that $ q_{ij}= a_{ji}$ for $i\ne j$  and $ \dis q_{jj}= -\sum_{k=1}^{\ell}a_{kj}.$	Hence the system $ \dis \frac{d\,\mbn(t)}{dt}=\nu D\mbn(t)$ is equivalent to $ \dis \frac{d\,\mbn^{\mtt}(t)}{dt}= \nu\mbn^{\mtt}(t) Q$.
Notice that $Q$ is the infinitesimal generator of an irreducible  Markov process with state space  $\big\{\; \dfrac{n}{\ell}, n=0, 1,\cdots, \ell\; \big\}^{\ell} .$  By using Theorem 5.1 of Pardoux \cite{par}, it follows that there exists a unique strictly positive equilibrium $\mbn^*$ which solves the equation $(\mbn^*)^{\!\mtt}Q=0 $.
	Since the state space is finite, the Markov process associed to the infinitesimal generator $Q$ is reccurent, and then the global asymptotic stability of $ \mbn^* $ is garanteed by the Theorem 6.5 of  the above reference.
	\fpr
\end{pr}

Next we treat the  existence and stablity of the endemic equilibrium for the system given by (\ref{deta}). Notice that (\ref{deta}) is of triangular form, and hence the theory of asymptotically autonomous systems for triangular systems (Vidyasagar \cite{vi}\,) guarantees that the asymptotic stability of its equilibrium is equivalent to that of the system
$$ \dfrac{d\,\mi(t)}{dt}= \nu\md\,\mi(t)- \ma\,\mi(t)+ \Big(\mbi_\ell-\text{diag}(1/\mbn_j^*)\text{diag}(\mi(t))\Big)\mb\,\mi(t),$$ 
where $ ( \mbn_1^*, \cdots , \mbn_{\ell}^* ) $ is the unique asymptotically stable equilibrium defined in \ref{stn}.
The Jacobian matrices  of $\Big(\mbi_\ell-\text{diag}(\mbn_j^{-1}(t))\text{diag}(\mi(t))\Big)\mb\,\mi(t)$  and $\nu \md\,\mi(t)- \ma\,\mi(t) $, respectively, at the  disease free equilibrium are $\mb$  and $ \mv$. For the convenience of the reader, we recall the following result.

\begin{thm}[Vidyasagar\cite{vi},\ Theorem 3.1  and  3.4]\label{syr} ~
	
Let $f$ and $g$ be two functions of class $\mathcal{C}^1$.
Consider  the following  system
	\begin{equation}\label{vthm}
	\hspace{-3cm}\left\{ 
	\begin{aligned} 
	\dot{x} &= f(x) \\
	\dot{y} & = g(x,y) \hspace{3cm} x\in \mathbb{R}^n, y\in\mathbb{R}^m\\
	&\text{with an equilibrium point } \, (x^*,y^*), \; \;  \text{i.e.}\, ,\\
	&f(x^*)=0 \; \text{and}\;  g(x^*,y^*)=0.
	\end{aligned} 
	\right.
	\end{equation}
	If $x^*$ is globally asymptotically stable (GAS) in $\mathbb{R}^n$ for the system $ \dot{x}= f(x)$, and if $y^*$  is GAS in $\mathbb{R}^m$ for  the system $ \dot{y}= g(x^*,y)$, then $(x^*,y^*)$ is locally asymptotically stable for (\ref{vthm}). Moreover, if all the trajectories of (\ref{vthm}) are forward  bounded, then $(x^*,y^*)$ is GAS for  (\ref{vthm}).
\end{thm}
We shall need below the
\begin{thm}[Hirsch \cite{h84}, Theorem 6.1]\label{h}~\\
	Let F be a $ \mathcal{C}^1$  vector field in  $ \mathbb{R}^q$,  whose flow $\phi$ preserves $ \mathbb{R}^q_+$ for $t\ge 0$ and is strongly monotone in $ \mathbb{R}^q_+$. Assume that the origin is an equilibrium and that all trajectories in $ \mathbb{R}^q_+$ are bounded. If the matrix-valued map
	$ \mathcal{D}F : \mathbb{R}^q\longrightarrow \mathbb{R}^q\times  \mathbb{R}^q$ is strictly decreasing, in the sense that
	$$ \text{if} \; \; x < y \; \; \text{then} \; \; \mathcal{D}F(x)> \mathcal{D}F(y), $$
	then either all trajectories in $ \mathbb{R}^q_+\backslash \{0\} $ tend to the origin, or there is a unique equilibrium $p^*$, ($p^*\gg 0$) in the interior of $ \mathbb{R}^q_+$ and all trajectories in $ \mathbb{R}^q_+\backslash\{0\}$  tend to $p^*$.
\end{thm}
Now, we are in a position to prove the main result of this section.
\begin{thm}\label{see}$ \bf{[Existence \ and \ stability \ of \ the \ EE]}$~\\
	Assume  that $\nu_{I}=\nu_{S}:=\nu$ and $\mathcal{R}_0 >1$. Then the system (\ref{md}) has a unique endemic equilibrium $z^*=\Big(\ms_1^* ,  \mi_1^*, \ms_2^* , \mi_2^*\cdots, \ms_{\ell}^*, \mi_{\ell}^*   \Big),$  which is globally asymptotically stable.
\end{thm}

\begin{pr} 
	It follows from  \ref{syr} that it is sufficient to study the stability of the reduced system
	$$ \dfrac{d\,\mi(t)}{dt} = \nu \md\,\mi(t)- \ma\,\mi(t)+ \Big(\mbi_\ell-\text{diag}(1/\mbn_j^*)\text{diag}(\mi(t))\Big)\mb\,\mi(t).$$ 
	Note that the set defined by 
	$$ K = \big\{ \big((u_1, \cdots, u_{\ell}),(v_1, \cdots, v_{\ell})\big)\in\mathbb{R}_+^{\ell}\times\mathbb{R}_+^{\ell}: 0\le v_i\le u_i, 1 \le  i  \le  \ell \; \text{and}\; \sum_{i=1}^{\ell}\!u_i=1\big\} $$
	is a compact  positively invariant for the system (\ref{detn}).
	Define  $$\mathcal{L}(\mi)= (\mb+\mv )\mi -\text{diag}(1/\mbn_j^*)\text{diag}(\mi)\mb\,\mi .$$ The derivative $ \mathcal{D}\mathcal{L}(\mi)$~is 
	\begin{eqnarray}  
	\mathcal{D}\mathcal{L}(\mi) &=& (\mb+\mv ) - \text{diag}(1/\mbn_j^*)\text{diag}(\mi)\mb- \text{diag}(1/\mbn_j^*)\text{diag}(\mb\,\mi) \nonumber\\
	&=& \nu \md-\ma +\mb - \text{diag}(1/\mbn_j^*)\text{diag}(\mi)\mb- \text{diag}(1/\mbn_j^*)\text{diag}(\mb\,\mi). \nonumber
	\end{eqnarray} 
	Notice that $\mathcal{D}\mathcal{L}(\rmi)$ is an irreducible Metzler matrix. Since $ \mb \ge 0$ and $\mb\ne 0 $, $\mathcal{D}\mathcal{L}$ is clearly stricly decreasing with respect of $\mi$. Applying  \ref{h}, we deduce that either all trajectories in K tend to the origin, or there is a unique equlibrium in the interior of K and all trajectories in $K\backslash([0,\infty)^{\ell}\times \{0\}^{\ell})$ tend to this equilibrium.\\
	We introduce the stability modulus $ \alpha(M)$ of a matrix $M$, which is the largest real part of the elements of the spectrum Spec(M) of M :
	$$ \alpha(M)= \un{\delta\in \text{Spec}(M)}{\max} \text{Re}(\delta).$$
	
	From Theorem 3.13 of Varga \cite{var} (chapter 3), $\mathcal{R}_0 > 1$ is equivalent to $ \alpha ( \mb+\mv)>0$, and the disease free equilibrium is unstable in this case. It then follows from \ref{h} that there exists a unique attracting endemic equilibrium $ \mi^*\neq 0 $, which satisfies
	\begin{eqnarray}\label{eee}
	(\nu\md-\ma+\mb)\mi^* -\text{diag}(1/\mbn_j^*)\text{diag}(\mi^*)\mb \,\mi^* &=& 0 .
	\end{eqnarray}
	Since $\mb$ is a non-negative matrix and $\mi^*\neq 0,$ by using (\ref{eee}), it follows that 
	\begin{eqnarray}\label{eqe}
	\mathcal{D}\mathcal{L}(\rmi^*)\mi^*= -\text{diag}(1/\mbn_j^*)\text{diag}(\mb\mi^*) \,\mi^* < 0 .
	\end{eqnarray}
	Using the fact that $ \mathcal{D}\mathcal{L}(\mi^*)$ is a Metzler matrix, (\ref{eqe}) implies that it is stable (Berman $\&$ Plemmons \cite{ber}: criteron $I_{28}$ of Theorem 6.2.3). The stability modulus then satisfies $\alpha\big(\mathcal{D}\mathcal{L}(\mi^*)\big)<0$. This proves the local asymptotic stability of $ \mi^*$. Since the attractivity of $\mi^*$ is guaranteed by Hirsh’s \ref{h}, we conclude that the endemic equilibrium $\mi^*$ is globally asymptotically stable if $\mathcal{R}_0> 1$.
	\fpr
\end{pr}
Let us mention that, under the assumption $ \nu_{I}=\nu_{S}:=\nu$, the DFE is globally assymptotically stable when $ \mathcal{R}_0= 1$. Indeed, $ \mathcal{R}_0= 1$ is equivalent to $ \alpha(\mb+\mv)=0$. But since $ \mb+\mv$ is an irreducible Metzler matrix, there exists a positive vector $v$ such that $(\mb+\mv)^{\mtt}v=0.$ Let us consider the Lyapunov function $ \dis \mathtt{L}(\mi)=\langle\; v\,|\,\mi\;\rangle $. The derivative of this function is 
\begin{eqnarray}
	\dot{\mathtt{L}}(\mi)&=& \langle\;  v\,|\,\dot{\mi}\; \rangle\n\\
	&=& \langle\; v\,|\,\mb+\mv- \text{diag}(1/\mbn_j^*)\text{diag}(\mi(t))\mb\,\mi(t)\; \rangle\n\\
	&=& -\langle\; v\,|\, \text{diag}(1/\mbn_j^*)\text{diag}(\mi(t))\mb\,\mi(t)\; \rangle\n\\
	&\le& 0. \n
\end{eqnarray}
Moreover, $ \dot{\mathtt{L}}(\mi)=0$ only at the DFE. Hence the DFE is GAS if $ \mathcal{R}_0= 1$. \\

\section{Comparison of the equilibra: connected patches model, homogeneous model, isolated patches}

We now consider the deterministic model in an homogeneous community:

 \begin{equation}\label{sh}
\hspace{-2cm}\left\{ 
\begin{aligned} 
\frac{d\,\ms}{dt}(t) & \!=\!  - \lambda\ms(t)\mi(t) +\gamma\,\mi(t) 
\\[2mm] 
\frac{d\,\mi}{dt}(t) & \! =\!  \lambda\ms(t)\mi(t) -\gamma\,\mi(t),
\end{aligned} \right.
\end{equation}
where $\lambda$ (resp. $\gamma$) is  the rate of the disease transmission (resp. recovery).\\ The endemic equilibrium of the system of ODEs (\ref{sh}) is $\mathfrak{z}^*=\left(\dfrac{\gamma}{\lambda}\, , 1-\dfrac{\gamma}{\lambda}\right)$.

We wish to compare this EE with  the one of the deterministic heterogeneous model.

\smallskip

$\bullet$ First, if the disease transmission and recovery rates are the same on all patches, that is for all $ j=1,\cdots, \ell$,  $\lambda_j=\lambda$ and $\gamma_j=\gamma$, then the EE of the ODEs (\ref{md}) is $ z^*=\left(\dfrac{\gamma}{\ell\lambda}, \dfrac{1}{\ell}\big(1-\dfrac{\gamma}{\lambda}\big),\cdots,\dfrac{\gamma}{\ell\lambda}, \dfrac{1}{\ell}\big(1-\dfrac{\gamma}{\lambda}\big) \right)$. In this case, we note  that  the  proportion of the infectious subpopulation in the hemogeneous model is equally distributed between all patches.

\smallskip

$\bullet$ We now look at the case where the patches have different disease transmission and recovery rates. In this case it is difficult to obtain the EE, even for a small number of patches. But it can be found relatively simply using any numerical solver when the state space is small. Here, we consider the case of two patches and use the solver "Wolfram Alpha" to compute the EE.

In the  below table we give the values of the infectious subpopulation in each patch at the equilibrium for several values  of the  parameters. We take $\gamma_1=\gamma_2=1$ and consider three cases. First $\lambda_1= 1.5$, $\lambda_2=2$. In this case when the patches are isolated, the value of the infectious subpopulation in patch 1 and patch 2 are $ \mi_1^*\approx 0.333$ and $\mi_2^*\approx 0.5$, respectively. Secondly $\lambda_1= 3$, $\lambda_2=2.5$, and $\mi_1^*\approx 0.666$ and $\mi_2^*\approx 0.600$ in isolated patches. Finally, in the case $\lambda_1= 1.5$, $\lambda_2=1.2$, we have $\mi_1^*\approx 0.333$ and $\mi_2^*\approx 0.166$.

\begin{table}[H]
	\begin{center}
	{\renewcommand{\arraystretch}{1.2} 
		{\setlength{\tabcolsep}{0.4cm} 
			\begin{tabular}{|c|c|c|c|c|c|c|}
				\hline  $\lambda_1$ & $\lambda_2$ & $\gamma_1$ & $\gamma_2$ & $\nu_I$ & $\nu_S$ & $\left(\dfrac{\mi_1^*}{\ms_1^*+\mi_1^*},\dfrac{\mi_2^*}{\ms_2^*+\mi_2^*}\right)$  \\
				\hline \rowcolor{gray!40} 1.5&2&1&1&0.0001&0.0001&  (0.332 , 0.507) \\
				\hline \rowcolor{gray!40} 1.5&2&1&1&0.0001&0.0005&  (0.334, 0.497) \\
				\hline \rowcolor{gray!40} 1.5&2&1&1&0.001&0.0001&  (0.333, 0.497) \\
				\hline \rowcolor{gray!40} 1.5&2&1&1&0.0001&0.001&  (0.332, 0.497)\\
				\hline \rowcolor{gray!100} &&&&&&  \\
				\hline \rowcolor{gray!40} 3&2.5&1&1&0.0001&0.0001&(0.667 , 0.598)   \\
				\hline \rowcolor{gray!40} 3&2.5&1&1&0.0007&0.0001&(0.666 , 0.599) \\
				\hline \rowcolor{gray!40} 3&2.5&1&1&0.001&0.0001&(0.666 , 0.598) \\
				\hline \rowcolor{gray!40} 3&2.5&1&1&0.0001&0.001&(0.666 , 0.598) \\
				\hline \rowcolor{gray!100} &&&&&&  \\
				\hline \rowcolor{gray!40} 1.5&1.2&01&1&0.0001&0.0001&(0.332 , 0.165)  \\
				\hline \rowcolor{gray!40} 1.5&1.2&1&1&0.0001&0.0009&(0.332 , 0.165)\\
				\hline \rowcolor{gray!40} 1.5&1.2&1&1&0.001&0.0001&(0.333 , 0.165)\\
				\hline \rowcolor{gray!40} 1.5&1.2&1&1&0.0001&0.008&(0.332 , 0.165)\\
				
				\hline
	\end{tabular}}}
\caption{proportion of $\mi_1^*$ and $\mi_2^*$ when patches are connected}
\end{center}
\end{table}

 In the  above table, we have the proportions of the infectious subpopulation in each patch at the equilibrium, for some values of the diffusion coefficients. We observe that those proportions are very close to those when patches are isolated.

\bigskip

We have shown that the stochastic model is well approximated by a deterministic patch model. If  $\mathcal{R}_0>1$, the system of ODEs (\ref{md}) has  a unique endemic equilibrium which is globally asymptotically stable. Moreover considering the case of two patches, it appears that in the heterogeneity case,  the final size of the epidemic in each patch is close to that of isolated patches. \\
In a future  work, we will study the fluctuations of the stochastic model around its deterministic law of large numbers limit. \\

\bigskip

\textbf{Acknowledgements.}
This paper is part of the author’s Ph.D. thesis written under the supervision of Prof.
Etienne Pardoux, whom the author would like to thank for careful reading and plentiful suggestions that greatly improved the paper. 

\bigskip

\textbf{Funding}: This research was partially supported by a thesis scholarship from the government of Ivory Coast, and a
salary as instructor at University of Aix–Marseille.

\medskip

\textbf{Conflict of Interest}: The author declare that there is no conflict of interest.

\end{document}